\documentclass[a4paper,10pt]{article}
\usepackage{amssymb}
\usepackage{latexsym}
\usepackage{amsmath}
\usepackage{amsthm}
\usepackage{amsfonts}
\usepackage{amssymb}
\usepackage{amsmath,amscd}
\usepackage{graphicx}

\newtheorem{Th}{Theorem}
\newtheorem{Prop}{Proposition}

\newtheorem{Con}{Conjecture}
\newcommand{\be}{\begin{equation}}
\newcommand{\ee}{\end{equation}}
\newcommand{\bes}{\begin{equation*}}
\newcommand{\ees}{\end{equation*}}

\newcommand{\R}{\mathbb{R}}

\newcommand{\C}{\mathbb{C}}
\newcommand{\Z}{\mathbb{Z}}

\newcommand\res{\mathop{\hbox{\vrule height 7pt width .5pt depth 0pt
\vrule height .5pt width 6pt depth 0pt}}\nolimits}

\newcommand{\reset}{\setcounter{equation}{0}\setcounter{Th}{0}\setcounter{Prop}{0}\setcounter{Co}{0}
\setcounter{Lm}{0}\setcounter{Rm}{0}}

\def\ti{\tilde}
\def\lf{\left}
\def\rg{\right}

\def\al{\alpha}
\def\la{\lambda}

\def\ep{\varepsilon}

\def\ds{\displaystyle}
\def\ov{\overline}

\def\p{\partial}

\def\res{\mathop{\hbox{\vrule height 7pt width .5pt 
depth 0pt\vrule height .5pt width 6pt depth 0pt}}\nolimits}

\DeclareMathOperator{\dist}{dist}

\begin{document}

\title{ Harmonic Maps from $S^3$ into $S^2$ with low Morse Index}

\author{ Tristan Rivi\`ere\footnote{Department of Mathematics, ETH Zentrum,
CH-8093 Z\"urich, Switzerland.}}

%\date{ }
\maketitle

{\bf Abstract :}{\it We prove that any smooth harmonic map from $S^3$ into $S^2$ of Morse index less or equal than 4 has to be an harmonic morphism, that is the successive composition of an isometry of $S^3$, the Hopf fibration and an holomorphic map from ${\C}P^1$ into itself.}

\medskip

\noindent{\bf Math. Class. 53C43, 58E20, 58J05, 35A15, 35J20  }

\section{Introduction}

Rigidity phenomena under Morse index control is an important problematic of the calculus of variations with applications in differential geometry. The archetype result
illustrating this question is the rigidity theorem of Francisco Urbano asserting that the closed minimal surfaces of Morse index less or equal than 5 in the canonical 3-sphere $S^3$ are given exclusively by geodesic spheres and Clifford torii.

In the present work we are considering smooth harmonic maps from the canonical $3-$sphere $S^3$ of ${\R}^4$ into the canonical $2-$sphere $S^2$ of ${\R}^3$ that is $C^\infty$ maps, critical point of
\[
E(u):=\frac{1}{2}\int_{S^3}|du|^2\ dvol_{S^3}
\]
in $W^{1,2}(S^3,S^2)$. The second derivative of $E$ at a critical point $u$ is given by
\[
D^2E_u(w):=\frac{1}{2}\int_{S^3}|dw|^2-|w|^2\,|du|^2\ dvol_{S^3}
\]
for $w$ a section of the bundle $u^{-1}TS^2$. The Morse index of $D^2E$ is the
dimension of the space spanned by the  eigen-sections with negative eigenvalues  in  the space of sections of $u^{-1}TS^2$ denoted  $\Gamma(u^{-1}TS^2)$ of
\[
L_u(w):=-P_u\Delta_{S^3} w-w\, |du|^2
\]
where $P_u$ is the projection map onto $T_uS^2$. Our main result in the present work is the following.
\begin{Th}
\label{th-I.1}
Let $u$ be a smooth harmonic map from $S^3$ into $S^2$. Assume that the Index of the quadratic form $D^2E_u$ on $\Gamma(u^{-1}TS^2)$ is less or equal than $4$ then it is equal to $4$ and there exists an isometry $\Phi$ of $S^3$ and an holomorphic map from ${\C}P^1$ into ${\C}P^1$ such that
\[
u=v\circ {\frak h}\circ \Phi\quad,
\]
where $\frak h$ is the Hopf map given by ${\frak h}(z,w)=(2\,z\,\ov{w},|z|^2-|w|^2)$.
\hfill $\Box$
\end{Th}
Since the work of Ahmed El Soufi it is known that the Morse index of any \underbar{smooth} harmonic map from $S^3$ into $S^2$ is larger or equal than 4. We conjecture that this result does not extend to weakly harmonic map and 4 should be replaced by 3 in this more general framework (see some comments in section IV). It is proved in \cite{Ura} that the index of the Hopf map is 4.

The theorem~\ref{th-I.1} has been elaborated in the course of the development of the following research program which is exposed in more details in the Master Thesis of Yujie Wu \cite{Wu}. We introduce the non empty subfamily of $W^{1,2}(S^3,S^2)$ given by
\[
{\mathcal A}:=\lf\{
\begin{array}{c}
\ds u\in C^0(B^4,W^{1,2}(S^3,S^2))\quad;\quad\max_{a\in \p B^4} \|du_a\|_{L^2(S^3)}\le\delta\\[5mm]
\ds\quad\mbox{ and }a\in\p B^4\rightarrow\frac{\ds\int_{S^3}u_a\ dvol_{S^3}}{\ds\lf| \int_{S^3}u_a\ dvol_{S^3}\rg|}\in S^2\mbox{ is non zero homotopic}
\end{array}
\rg\}
\]
where $\delta>0$ is chosen small enough in such a way that, using Poincar\'e inequality, 
\[
\|du_a\|_{L^2(S^3)}\le\delta\quad\Longrightarrow\quad {\ds\lf| \int_{S^3}u_a\ dvol_{S^3}\rg|}>1/2
\]
It is proved in \cite{Wu} that  for $\delta>0$ again chosen small enough
\be
\label{minmax}
W_{\mathcal A}:=\inf_{u\in {\mathcal A}}\max_{a\in B^4}E(u_a)>\delta\quad.
\ee
We conjecture that the minmax is realized by harmonic morphisms and precisely
\[
W_{\mathcal A}=\frac{1}{2}\int_{S^3}|d\frak h|^2\ dvol_{S^3}\quad.
\]
If this is true then  a proof of the following easily follows.
\begin{Con}\cite{Riv}
\label{conj-riv}
The minimizers of the 3-energy among homotopically non zero maps from $S^3$ into $S^2$ are exactly given by the composition
of the Hopf maps $\frak h$ with conformal transformations and isometries of $S^3$ .
\end{Con}
The overall framework of this program shares a lot of similarities with the structure of the proof of the Willmore Conjecture by Fernando Cod\'a Marques and Andr\'e Neves \cite{MN}.
The counterparts of the area, the Willmore energy, the canonical family and Urbano's result are respectively given in our program by the 2-energy, the 3-energy, the family $u\circ\phi_a$ (see \ref{moebius}) and the theorem~\ref{th-I.1}.

The framework of harmonic maps offers interesting new phenomena which are not present in the corresponding questions for minimal surfaces. This is related to the fact that not every harmonic map from $S^3$ into $S^2$ is strongly approximable by smooth maps. This is mentioned in section IV.
\section{Preliminaries.}
\reset
\subsection{Global Frame on $S^3$.}

We consider the following positive orthonormal frame of $S^3$ which extends as a free orthogonal family in the whole ${\R}^4$ :
\[
\lf\{
\begin{array}{l}
\ds e_1:=x_1\, \p_{x_2}-x_2\, \p_{x_1}+x_3\, \p_{x_4}-x_4\, \p_{x_3}\\[5mm]
\ds e_2:=x_2\, \p_{x_3}-x_3\, \p_{x_2}+x_1\, \p_{x_4}-x_4\, \p_{x_1}\\[5mm]
\ds e_3:=x_3\, \p_{x_1}-x_1\, \p_{x_3}+x_2\, \p_{x_4}-x_4\, \p_{x_2}
\end{array}
\rg.
\]
In other words using cyclic indexation in ${\Z}_3$ we have
\be
\label{II.1}
e_i:= x_i\,\p_{x_{i+1}}-x_{i+1}\,\p_{x_i}+x_{i-1}\, \p_{x_4}-x_4\, \p_{x_{i-1}}
\ee
This gives in particular
\be
\label{II.1-a}
dx_i=-x_{i+1}\ e_i^\ast-x_4\, e^\ast_{i+1}+x_{i-1}\ e^\ast_{i-1}\quad\mbox{ and }\quad dx_4=\sum_{i=1}^4 x_{i-1}\ e_i^\ast
\ee
We have respectively
\be
\label{II.2}
[e_i,e_{i+1}]:=\sum_{k,l=1}^4 \lf(e_i^k\, \p_{x_k}e_{i+1}^l-e_{i+1}^k\, \p_{x_k}e_i^l\rg)\ \p_{x_l}=-2\, e_{i-1}
\ee
We denote
\[
e^\ast_i:=\iota_{\p B^4}^\ast \lf(x_i\,d{x_{i+1}}-x_{i+1}\,d{x_i}+x_{i-1}\, d{x_4}-x_4\, d{x_{i-1}}\rg)
\]
We have
\be
\label{II.3}
\begin{array}{l}
<e^\ast_i, e_j>=<x_i\,d{x_{i+1}}-x_{i+1}\,d{x_i}+x_{i-1}\, d{x_4}-x_4\, d{x_{i-1}}, (\iota_{\p B^4})_\ast e_j>\\[5mm]
\quad= \lf<x_i\,d{x_{i+1}}-x_{i+1}\,d{x_i}+x_{i-1}\, d{x_4}-x_4\, d{x_{i-1}}, x_j\,\p_{x_{j+1}}-x_{j+1}\,\p_{x_j}+x_{j-1}\, \p_{x_4}-x_4\, \p_{x_{j-1}}\rg>=\delta_{ij}
\end{array}
\ee
Hence $(e_i^\ast)_{i=1,2,3}$ realizes a dual frame to $(e_i)_{i=1,2,3}$. We have in one hand
\be
\label{II.4}
de_i^\ast=2\,\iota_{\p B^4}^\ast \lf(dx_i\wedge dx_{i+1}+dx_{i-1}\wedge dx_4\rg)
\ee
and in the other hand, explicit computations (see \cite{Riv}) give
\be
\label{II.5}
e^\ast_{i+1}\wedge e^\ast_{i-1}=\iota_{\p B^4}^\ast \lf(dx_i\wedge dx_{i+1}+dx_{i-1}\wedge dx_4\rg)
\ee
Hence we deduce
\be
\label{II.6}
de_i^\ast=2\, e^\ast_{i+1}\wedge e^\ast_{i-1}
\ee
Denote by $\nabla$ the Levi Civita covariant derivative on $S^3$. Since $\nabla$ is torsion free we have
\[
\nabla_{e_i}e_{i+1}=\nabla_{e_{i+1}}e_i+[e_i,e_{i+1}]=\nabla_{e_{i+1}}e_i-2\, e_{i-1}\quad.
\]
We have also
\[
(\nabla_{e_i}e_{i+1}, e_{i+1})\equiv 0\quad\mbox{ and }\quad(\nabla_{e_{i+1}}e_{i}, e_{i})\equiv 0\quad.
\]
Hence
\[
\nabla_{e_i}e_{i+1}=\la_i\, e_{i-1}\quad\mbox{ and }\quad \nabla_{e_{i+1}}e_{i}=(\la_i+2)\, e_{i-1}\quad.
\]
Since
\[
(\nabla_{e_i}e_{i+1}, e_{i-1})+(\nabla_{e_i}e_{i-1}, e_{i+1})=0\quad.
\]
We deduce that
\[
\la_i+\la_{i-1}+2=0\quad\forall i\in {\Z}_3\quad.
\]
This implies that $\la_i=-1$ and we deduce finally
\be
\label{II.7}
\nabla_{e_i}e_{i+1}=-\, e_{i-1}\quad\mbox{ and }\quad \nabla_{e_{i+1}}e_{i}=\, e_{i-1}\quad.
\ee
Since
\[
(\nabla_{e_i}e_i,e_{i+1})+(e_i,\nabla_{e_i}e_{i+1})=0\quad\mbox{ and }\quad(\nabla_{e_i}e_i,e_{i-1})+(e_i,\nabla_{e_i}e_{i-1})=0\quad,
\]
we deduce from the previous that
\be
\label{II.8}
\nabla_{e_i}e_i\equiv 0\quad.
\ee
Let $\al$ be a one form, we have for any pair of vector field $X$ and $Y$
\be
\label{II.8a}
\lf< \nabla_X\al,Y \rg>=-\lf<\al,\nabla_XY\rg>+d_X\lf(<\al,Y>\rg)
\ee
From this identity we deduce respectively
\be
\label{II8b}
\nabla_{e_i}e_i^\ast=0\quad,\quad\nabla_{e_{i+1}}e_i^\ast=e_{i-1}^\ast\quad\mbox{ and }\quad \nabla_{e_{i-1}}e_i^\ast=-e_{i+1}^\ast
\ee
\subsection{Conformal Killing Fields}
We denote for $i=1\cdots 4$ by $X^j$ the following conformal killing vector field on $S^3$
\be
\label{II8c}
X^j(x):=\ep_j-x_j\, x\quad,
\ee
where $(\ep_j)_{j=1\cdots 4}$ is the canonical basis of ${\R}^4$. We consider in particular $X^1$ that we express in the $e$-frame. Since $e_i(x)$ is orthogonal to $x$ we have
\[
X^1=\sum_{i=1}^3<\ep_1, e_i>\ e_i=-x_2\, e_1-x_4\, e_2+x_3\, e_3\quad.
\]
Observe that
\[
X^1=(dx_1)^\ast\quad.
\]
Similarly we could have considered $X^i:= (dx_i)^\ast$.
Hence, using (\ref{II.7}) we obtain successively
\be
\label{II.9}
\begin{array}{l}
\ds\nabla_{e_1}X^1=-<dx_2,e_1>\, e_1-<dx_4,e_1>\, e_2+<dx_3,e_1>\, e_3-x_4\ \nabla_{e_1}e_2+x_3\, \nabla_{e_1}e_3\\[5mm]
\ds\quad=-x_1\, e_1-x_3\, e_2-x_4\, e_3+x_4 \,e_3+x_3\, e_2=-x_1\, e_1\quad,
\end{array}
\ee
then
\be
\label{II.10}
\begin{array}{l}
\ds\nabla_{e_2}X^1=-<dx_2,e_2>\, e_1-<dx_4,e_2>\, e_2+<dx_3,e_2>\, e_3-x_2\ \nabla_{e_2}e_1+x_3\, \nabla_{e_2}e_3\\[5mm]
\ds\quad=x_3\, e_1-x_1\, e_2+x_2\, e_3-x_2 \,e_3-x_3\, e_1=-x_1\, e_2\quad,
\end{array}
\ee
and
\be
\label{II.11}
\begin{array}{l}
\ds\nabla_{e_3}X^1=-<dx_2,e_3>\, e_1-<dx_4,e_3>\, e_2+<dx_3,e_3>\, e_3-x_2\ \nabla_{e_3}e_1-x_4\ \nabla_{e_3}e_2\\[5mm]
\ds\quad=x_4\, e_1-x_2\, e_2-x_1\, e_3+x_2 \,e_2-x_4\, e_1=-x_1\, e_3\quad.
\end{array}
\ee
Combining these facts one obtains in particular
\be
\label{II.12}
\nabla_{e_i}X^1=-x_1\, e_i\quad\mbox{ and }\quad \sum_{i=1}^3\nabla_{e_i}\nabla_{e_i}X^1=-\sum_{i=1}^3 <dx_1, e_i> \ e_i= x_2\, e_1+x_4\ e_2- x_3\ e_3=-\, X^1
\ee
\subsection{Weitzenb\"ock Formula}
Under the notations above one has the following expression of the exterior differential of a $p-$form $\al$ on $S^3$
\be
\label{d-charac-nabla}
d\al=\sum_{i=1}^3 e_i^\ast\wedge\nabla_{e_i}\al\quad.
\ee
We recall the definition of the interior product of a vector field $X$ with a $p-$form $\al$
\[
(X\res \al)(Y_1\cdots Y_{p-1}):=\al(X,Y_1\cdots Y_{p-1})\quad.
\]
Finally we recall the definition of the Ricci tensor acting on multi-vectors or differential forms
\[
R(X,Y)\,Z:=\nabla_X\nabla_YZ-\nabla_Y\nabla_XZ-\nabla_{[X,Y]}Z
\]
We have in particular
\be
\label{II.13}
R(e_i,e_{i+1})\,e_i=-e_{i+1}\quad,\quad R(e_i,e_{i+1})\,e_{i+1}=e_i\quad\mbox{ and }\quad R(e_i,e_{i+1})\,e_{i-1}=0\quad.
\ee
and
\be
\label{II.14}
R(e_i,e_{i+1})\,e^\ast_i=-e^\ast_{i+1}\quad,\quad R(e_i,e_{i+1})\,e^\ast_{i+1}=e^\ast_i\quad\mbox{ and }\quad R(e_i,e_{i+1})\,e_{i-1}^\ast=0\quad.
\ee
The Weitzenb\"ock formula reads as follows. Let $\al$ be a $p-$form, we have
\be
\label{II.15}
(dd^\ast+d^\ast d)\al=-\sum_{i=1}^3\nabla_{e_i}\nabla_{e_i}\al-\sum_{ij=1}^3e_i^\ast\wedge\lf(e_j\res\lf(R(e_i,e_j)\,\al\rg)\rg)
\ee
\subsection{The second Derivative of the Dirichlet Energy}
Let $u$ be a smooth map from $S^3$ into $S^2$ and let $w$ be a smooth map from $S^3$ into ${\R}^3$ whose $L^\infty$ norm is less than $1/4$ with also small $C^1$ norm. We compute
\[
\begin{array}{l}
\ds E\lf( \frac{u+w}{|u+w|} \rg)=\frac{1}{2}\int_{S^3}\lf|d \frac{u+w}{|u+w|}\rg|^2\ dvol_{S^3}\\[5mm]
\ds \quad=E(u)+\int_{S^3}du\cdot dw-u\cdot w\,|du|^2\ dvol_{S^3}+\frac{1}{2}\int_{S^3}|dw|^2-|w|^2\,|du|^2\ dvol_{S^3}\\[5mm]
\ds\quad+\frac{1}{2}\int_{S^3}\lf[3\, (u\cdot w)^2\, |du|^2-\, 4\, du\cdot dw\, (u\cdot w)-2\, |d(u\cdot w)|^2\rg]\ \ dvol_{S^3}+O((|w|+|dw|)^3)
\end{array}
\]
Assuming $u$ is harmonic and $w$ is a section of $u^{-1}TS^2$ that is to say $u\cdot w\equiv 0$ we then have
\be
\label{II.16}
E\lf( \frac{u+w}{|u+w|} \rg)=E(u)+\frac{1}{2}\int_{S^3}|dw|^2-|w|^2\,|du|^2\ dvol_{S^3}+O((|w|+|dw|)^3)\quad.
\ee
Let $P_u$ be the map from $S^3$ into $3$ by $3$ matrices giving at each $x$ the orthogonal projection onto $T_{u(x)}S^2$, we are interested in the operator $L$
defined from $\Gamma(u^{-1}TS^2)$ into itself and given by
\[
L_u(w):=-P_u\Delta_{S^3} w-w\, |du|^2\quad,
\]
where $\Delta_{S^3}$ is the negative Laplace Beltrami operator on $S^3$. On $\Gamma(u^{-1}TS^2)$ the second derivative of $E$ is then given by
\[
D^2E_u(w)=\frac{1}{2}\int_{S^3}w\cdot L_u(w)\ dvol_{S^3}\quad.
\]
\subsection{The infinitesimal variation of $E$ along the action of conformal transformations of $S^3$}
Let $u$ be a smooth harmonic map from $S^3$ into $S^2$ and $X^i$ denote the conformal Killing fields given by (\ref{II8c}). We are interested in computing $L_u(<du,X^i>)$.
We have using respectively (\ref{II.6}) and (\ref{II.12})
\be
\label{II.17}
\begin{array}{l}
\ast d\ast d<du,X^1>=\ast d\ast <\nabla du,X^1>+\ast d\ast<du,\nabla X^1>\\[5mm]
\ds\quad=\ast\sum_{i=1}^3d\lf(e^\ast_{i+1}\wedge e_{i-1}^\ast\ <\nabla_{e_i} du,X^1>\rg)+\ast\sum_{i=1}^3d\lf(e^\ast_{i+1}\wedge e_{i-1}^\ast\ <du,\nabla_{e_i}X^1>\rg)\\[5mm]
\ds\quad=\sum_{i=1}^3<\nabla_{e_i}\nabla_{e_i} du,X^1>+2\, <\nabla_{e_i}du, \nabla_{e_i}X^1>+<du,\nabla_{e_i}\nabla_{e_i}X^1>\\[5mm]
\ds\quad=-<du,X^1>+\sum_{i=1}^3<\nabla_{e_i}\nabla_{e_i} du,X^1>-\,2\, x_1\,<\nabla_{e_i}du, e_i>
\end{array}
\ee
The Weitzenb\"ock Formula and the harmonicity of $u$ give
\be
\label{II.18}
\begin{array}{l}
\ds\sum_{i=1}^3\nabla_{e_i}\nabla_{e_i} du=- dd^\ast du-\sum_{ij=1}^3e_i^\ast\wedge\lf(e_j\res\lf(R(e_i,e_j)\,du\rg)\rg)\\[5mm]
\ds\quad=-du\, |du|^2-u\, d|du|^2-\sum_{ij=1}^3e_i^\ast\wedge\lf(e_j\res\lf(R(e_i,e_j)\,du\rg)\rg)
\end{array}
\ee
Using (\ref{II.14}) we compute
\be
\label{II.19}
\lf\{
\begin{array}{l}
e_i^\ast\wedge\lf(e_{i+1}\res\lf(R(e_i,e_{i+1})\,du\rg)\rg)=-\,e_i^\ast\, <du,e_i>\\[5mm]
\ds e_i^\ast\wedge\lf(e_{i-1}\res\lf(R(e_i,e_{i-1})\,du\rg)\rg)=-\,e_i^\ast\, <du,e_i>\quad.
\end{array}
\rg .
\ee
Hence combining (\ref{II.18}) and (\ref{II.19}) we obtain
\be
\label{II.20}
\sum_{i=1}^3\nabla_{e_i}\nabla_{e_i} du=-du\, |du|^2-u\, d|du|^2+\,2\, du\quad.
\ee
We have also
\be
\label{II.21}
\begin{array}{l}
\ds\sum_{i=1}^3<\nabla_{e_i}du, e_i>=\sum_{i=1}^3\lf<\nabla_{e_i}\lf(\nabla_{e_{i-1}}u\ e_{i-1}^\ast+\nabla_{e_{i}}u\ e_{i}^\ast+\nabla_{e_{i+1}}u\ e_{i+1}^\ast\rg),e_i\rg>\\[5mm]
\ds\quad=\sum_{i=1}^3\nabla_{e_i}\nabla_{e_i} u=\Delta_{S^3}u=-u\, |du|^2
\end{array}
\ee
Combining (\ref{II.17}), (\ref{II.20}) and (\ref{II.21}) we obtain for any $k=1\cdots 4$
\be
\label{II.22}
\begin{array}{l}
\Delta_{S^3}<du,X^k>=\ast d\ast d<du,X^k>= (-|du|^2+1)\, <du,X^k>+2\, x_k\, u\, |du|^2

\end{array}
\ee
This gives
\be
\label{II.23}
L(<du,X^1>)=-<du,X^1>
\ee
In particular we obtain for any $i=1\cdots 4$
\be
\label{II.24}
D^2E_u(<du,X^i>)=-\frac{1}{2}\int_{S^3}|<du,X^i>|^2\ dvol_{S^3}\quad.
\ee
Assume the $<du, X^i>$ are not linearly independent then there exist $X\in \mbox{Span}\{X^1\cdots X^4\}\setminus \{0\}$ such that $<du,X>\equiv 0$. Modulo the action of rotations we can assume $<du,X^1>\equiv 0$. This implies that the maps $u\circ\phi_{t\ep_1}$ are all constants where
\be
\label{moebius}
\phi_a(z):=(1-|a|^2)\frac{z-a}{|z-a|^2}-a\quad
\ee
As $t$  goes to $1$ the map $u\circ\phi_{t\ep_1}$ is converging in $W^{1,2}$ towards a constant map\footnote{Here we are using that the map $u$ is smooth, the singular harmonic morphism $u$ constructed in section IV does not satisfies this property.}. Hence $u$ has to be constant.
We have then given a proof of a result by Ahmed El Soufi\footnote{The result of El Soufi is a way more general and extends to general targets and arbitrary spheres
of dimensions strictly larger than 2.}.
\begin{Th}
\label{th-elsoufi} \cite{ElS}
Let $u$ be a smooth non constant harmonic map from $S^3$ into $S^2$ then the index of $u$ for the Dirichlet energy is at least $4$.\hfill $\Box$
\end{Th}
\subsection{Infinitesimal perturbations along the orthogonal of the push-forwards of Conformal Killing Fields.}

Recall the notation $(X^l)^\ast:=dx_l$ in $S^3$. Let $\al$ and $\beta$ two 1-forms on $S^3$ we claim that
\be
\label{II.41}
\al\cdot\beta=\sum_{l=1}^4 <\al,X^l>\, <\beta,X^l>
\ee
Indeed, denote $\pi(x):=x/|x|$. We have in particular $X^l=\pi_\ast\p_{x_l}$. We have
\[
\al\cdot\beta=\pi^\ast\al\cdot\pi^\ast\beta=\sum_{l=1}^4 <\pi^\ast\al,\p_{x_l}>\, <\pi^\ast\beta,\p_{x_l}>=\sum_{l=1}^4 <\al,\pi_\ast\p_{x_l}>\, <\beta,\pi_\ast\p_{x_l}>
\]
which gives (\ref{II.41}).

\medskip

We consider $w:=u\times<du,X^l>$ for $l=1\cdots 4$. We have using the harmonic map equation as well as (\ref{II.22}) and (\ref{II.41})
\[
\begin{array}{l}
\ds-P_u\Delta_{S^3}(u\times<du,X^l>)=-P_u(\Delta_{S^3}u\times<du,X^l>)-u\times\Delta_{S^3}<du,X^l>-2\,P_u(du\dot{\times}d<du,X^l>)\\[5mm]
\ds=u\times<du,X^l>\, |du|^2+(|du|^2-1)\,u\times<du,X^l>+2\, (du\times u)\cdot (du,<du,X^l>)\\[5mm]
\ds=(2\,|du|^2-1)\,u\times<du,X^l>-2\, \sum_{k=1}^4 u\times<du,X^k>\, <du,X^k>\cdot<du,X^l>
\end{array}
\]
Taking the scalar product with $u\times<du,X^l>$, summing over $l$ and integrating over $S^3$ gives
\be
\label{II.44}
\begin{array}{l}
\ds\int_{S^3}\sum_{l=1}^4|d(u\times<du,X^l>)|^2\ dvol_{S^3}-\int_{S^3}\sum_{l=1}^4|u\times<du,X^l>|^2\ [|du|^2-1]\ dvol_{S^3}\\[5mm]
\ds =\int_{S^3}\sum_{l=1}^4|u\times<du,X^l>|^2\ |du|^2\ dvol_{S^3}\\[5mm]
\ds-2\, \int_{S^3}\sum_{k,l=1}^4 (u\times<du,X^k>)\cdot (u\times<du,X^l>)\, <du,X^k>\cdot<du,X^l>\ dvol_{S^3}\\[5mm]
\ds=\int_{S^3}\sum_{k,l=1}^4 |<du,X^k>|^2\ |<du, X^l>|^2-2(<du,X^k>\cdot<du,X^l>)^2\ dvol_{S^3}\\[5mm]
\ds=\int_{S^3}-\sum_{k=1}^4 |<du,X^k>|^4\ +2\,\sum_{k< l}|<du,X^k>|^2\ |<du, X^l>|^2\ dvol_{S^3}\\[5mm]
\ds\quad\quad\quad\quad\quad\quad-4\int_{S^3}\sum_{k< l}(<du,X^k>\cdot<du,X^l>)^2\ dvol_{S^3}\\[5mm]
\ds=\int_{S^3}|du|^4-2\,|du\,\dot{\otimes}\,du|^2\ dvol_{S^3}
\end{array}
\ee
where in the 3 last identity we used (\ref{II.41}).
\section{Proof of Theorem~\ref{th-I.1}}
We are assuming that for any $w\in \Gamma(u^{-1}TS^2)$
\[
\int_{S^3}|dw|^2-[|du|^2-1]\, |w|^2\ dvol_{S^3}\ge 0\quad.
\]
Combining this assumption with (\ref{II.44}) gives in particular
\be
\label{II.45}
\int_{S^3}|du|^4-2\,|du\dot{\otimes}du|^2\ dvol_{S^3}\ge 0\quad.
\ee
At every point $du$ has at most rank 2. Assuming Rank$\,du_x=1$ we obviously get a contradiction. Let then $x$ be such that Rank$\,du_x=2$ and let $(f_1,f_2)$ be an orthonormal basis of the orthogonal 2-plane to Ker$\,du_x$. We have at $x$
\[
\begin{array}{l}
\ds|du|^4-2\,|du\dot{\otimes}du|^2=(|\p_{f_1}u|^2+|\p_{f_2}u|^2)^2-2\, |\p_{f_1} u|^4-2\, |\p_{f_2} u|^4-4\,|\p_{f_1}u\cdot\p_{f_2}u|^2\\[5mm]
\ds=-|\p_{f_1}u|^4-|\p_{f_2}u|^4+\, 2\, |\p_{f_1}u|^2\, |\p_{f_2}u|^2-4\,|\p_{f_1}u\cdot\p_{f_2}u|^2\\[5mm]
\ds=-(|\p_{f_1}u|^2-|\p_{f_2}u|^2)^2-4\,|\p_{f_1}u\cdot\p_{f_2}u|^2\le 0
\end{array}
\]
Hence we have proved that $|du|^4-2\,|du\dot{\otimes}du|^2\le 0$ at every point. The inequality (\ref{II.44}) implies then $|du|^4-2\,|du\dot{\otimes}du|^2\equiv 0$.
Hence at any point we have that, either $du_x=0$ or $du_x$ is transversally conformal. We conclude that $u$ defines an harmonic horizontally weakly conformal map
between $S^3$ and $S^2$ (see definition 2.4.2 in \cite{BW}). This gives that $u$ is an harmonic morphism between $S^3$ and $S^2$ (see lemma 4.2.1 from \cite{BW}) and using theorem 6.7.7 of \cite{BW} we deduce that $u$ is the successive composition of an isometry of $S^3$, the Hopf fibration and an holomorphic map from ${\C}P^1$ into itself. This concludes the proof of  theorem~\ref{th-I.1}.
\section{Non Smooth Harmonic Maps}
\reset
In this last section we would like to stress the importance between smooth and non smooth harmonic maps with respect to the above computations. 

\subsection{A non-balancing condition}

We will start by recalling a discrepancy regarding a  balancing condition discovered in \cite{St}.

\medskip

Taking the dot product of (\ref{II.22}) with $u$ and integrating over $S^3$ gives for $k=1\cdots 4$
\[
\begin{array}{l}
\ds2\,\int_{S^3}x_k\, |du|^2\ dvol_{S^3}=\int_{S^3} u\cdot\Delta_{S^3}<du,X^k>\ dvol_{S^3}\\[5mm]
\ds\quad=\int_{S^3}\Delta_{S^3}u\cdot <du,X^k>\ dvol_{S^3}=-\int_{S^3}u\cdot <du,X^k>\ |du|^2\ dvol_{S^3}=0
\end{array}
\]
We deduce then the following balancing condition first proved in \cite{St}
\begin{Prop}
\label{pr-balance}
Let $u$ be a smooth harmonic map from $S^3$ into $S^2$ then
\be
\label{balance}
\forall\, k=1\cdots 4\quad\quad\int_{S^3}x_k\ |du|^2_{S^3}\ dvol_{S^3}=0\quad.
\ee
\hfill $\Box$
\end{Prop}
Let $v$ be an harmonic map from $S^2$ into $S^2$. We extend radially $v$ into the map $\tilde{v}$ from ${\R}^3$ into $S^2$.
Recall the formula of the laplacian in polar coordinates where $\Delta_{S^2}$ is the negative Laplace Beltrami operator on the $2-$sphere
\[
\Delta_{{\R}^3}=\p^2_{r^2}+2\,r^{-1}\,\p_r+r^{-2}\Delta_{S^2}\quad.
\]
Since by definition $\p_r\ti{v}=0$ we have that $\Delta_{{\R}^3}\ti{v}=r^{-2}\Delta_{S^2}v$ and since $v$ is harmonic we have
\[
v\wedge\Delta_{S^2} v=0\quad\Rightarrow\quad \ti{v}\wedge\Delta_{{\R}^3}\ti{v}=0
\]
Hence $\ti{v}$ realizes an harmonic map from ${\R}^3$ into $S^2$ and satisfies the equation\footnote{Observe that the equation is also satisfied throughout the origin. Indeed, by classical theory $v$ is smooth then the following bound holds : $|x|\,|\nabla \ti{v}|=O(1)$ near the origin and then $\nabla \ti{v}\in L^{3,\infty}(B_1(0))$. The following equation holds obviously
\be
\label{harm}
\mbox{div}(\ti{v}\wedge\nabla\ti{v})=0\quad \mbox{in }{\mathcal D}'(B^3_1(0)\setminus\{0\})
\ee
Since $\ti{v}\wedge\nabla\ti{v}\in L^p(B^3_1(0))$ for some $p>3/2$, a classical capacity argument gives that (\ref{harm}) extends throughout the origin.}
\[
-\Delta_{{\R}^3}\ti{v}=\ti{v}\, |\nabla_{{\R}^3}\ti{v}|^2
\]
Denote $\ti{u}(x_1,x_2,x_3,x_4):=\ti{v}(x_1,x_2,x_3)=v(|x'|^{-1} x')$ where $x'=(x_1,x_2,x_3)$. We have for $\rho^2=x_1^2+x_2^2+x_3^2+x_4^2$
\[
\p_\rho\ti{u}=0 \quad\quad\mbox{in } {\R}^4\setminus\{x'=0\}\quad.
\]
Hence, away from the $x_4$ axis $\{x'=0\}$ we have
\[
\ti{u}(x)\wedge\rho^{-2}\Delta_{S^3}\ti{u}(x)=\ti{u}(x)\wedge\Delta_{{\R}^4}\ti{u}(x)=\ti{v}(x')\wedge\Delta_{{\R}^3}\ti{v}(x')=0
\]
We deduce that the restriction $u$ of $\ti{u}$ to the $3-$sphere away from the north and the south pole $x_4=\pm 1$ is a smooth harmonic map. We have obviously
\[
|x'|\,|\nabla_{S^3} u|(x)=O(1)
\]
Moreover we have respectively $\dist_{S^3}(x,\mbox{North})=\arcsin |x'|$ for $x_4>0$ and $\dist_{S^3}(x,\mbox{South})=\arcsin |x'|$ for $x_4<0$ hence  $\nabla_{S^3} u\in L^{3,\infty}(S^3)$. The capacity argument mentioned in the previous note gives that $u$ is weakly harmonic on $S^3$.

Take now $v_a(x'):=\phi_a(x')$ where $a\in B_1^3(0)$ and $\phi_a$ is given by (\ref{moebius}). It is a well known fact that this map realizes a conformal harmonic diffeomorphism of $S^2$ we have moreover
\[
\lim_{a\rightarrow (1,0,0)}\phi_a= -1\quad\quad\mbox{ in } C^1_{loc}(S^2\setminus\{(-1,0,0)\})
\]
and
\be
\label{radon-conv}
|\nabla_{S^2}\phi_a|^2\ dvol_{S^2}\rightharpoonup 8\pi\, \delta_{(-1,0,0)}\quad\mbox{ in Radon measure as }a\rightarrow (1,0,0)
\ee
Take $u_a(x):=v_a(x')$. The map $u_a$ realizes a {\it singular harmonic morphism} with the north and the south poles as singularities. We deduce from (\ref{radon-conv}) that
\be
\label{contrebalance}
\lim_{a\rightarrow (1,0,0)}\int_{S^3} x_1\, \ |du_a|^2_{S^3}\ dvol_{S^3}<0
\ee
This implies the following proposition.
\begin{Prop}
\label{pr-contrebalance}
There exists a weak harmonic map  from $S^3$ into $S^2$ such that
\be
\label{contre}
\int_{S^3}x_1\ |du|^2_{S^3}\ dvol_{S^3}\ne 0\quad.
\ee
\hfill $\Box$
\end{Prop}

\medskip

To conclude this paper we would like to give some hint why the main theorem~\ref{th-I.1} should not hold for non smooth harmonic maps. Let
\[
{\mathcal B}:=\lf\{
\begin{array}{c}
\ds u\in C^0(B^3,W^{1,2}(S^3,S^2))\quad;\quad\max_{b\in \p B^3} \|du_b\|_{L^2(S^3)}\le\delta\\[5mm]
\ds\quad\mbox{ and }b\in\p B^3\rightarrow\frac{\ds\int_{S^3}u_b\ dvol_{S^3}}{\ds\lf| \int_{S^3}u_b\ dvol_{S^3}\rg|}\in S^2\quad\mbox{ is non zero homotopic}
\end{array}
\rg\}
\]
where $\delta>0$ is chosen small enough in such a way that, using Poincar\'e inequality, 
\[
\|du_b\|_{L^2(S^3)}\le\delta\quad\Longrightarrow\quad {\ds\lf| \int_{S^3}u_b\ dvol_{S^3}\rg|}>1/2
\]
We claim that ${\mathcal B}$ is not empty. We take the harmonic map $u$ we just constructed in the first part of this section and consider the family
\[
b\in B^3=B^4\cap\{x\in B^4\ ;\ x_4=0\}\quad \longrightarrow\quad u\circ\phi_{(1-\ep)\,b}
\]
For any $\delta>0$ there exists $\ep>0$ small enough such that $\max_{b\in \p B^3} \|d(u\circ\phi_{(1-\ep)\,b})\|_{L^2(S^3)}\le\delta$ it is moreover not difficult to check that
\[
b\in\p B^3\longrightarrow\frac{\ds\int_{S^3}u\circ\phi_{(1-\ep)\,b}\ dvol_{S^3}}{\ds\lf| \int_{S^3}u\circ\phi_{(1-\ep)\,b}\ dvol_{S^3}\rg|}\in S^2\quad\mbox{ is non zero homotopic}\quad.
\]
Introducing 
\be
\label{II.46}
W_{\mathcal B}:=\inf_{u\in {\mathcal B}}\max_{b\in B^3}E(u_b)
\ee
one easily show that $W_{\mathcal B}>\delta$. It is expected that this minmax is achieved by an harmonic map of index 3. The singular harmonic morphism $u$ is a natural candidate for realizing this minmax.

\end{document}